\documentclass{ifacconf}

\usepackage{graphicx}      
\usepackage{natbib}        
\usepackage{amsmath}
\usepackage{amsfonts}
\usepackage{xcolor}
\usepackage{booktabs}

\DeclareMathOperator{\trace}{Trace}

\begin{document}
\begin{frontmatter}

\title{A distributed voltage stability margin for power distribution networks} 

\thanks[footnoteinfo]{This research is supported by ETH funds and the SNF Assistant
Professor Energy Grant \# 160573.}

\author[First]{Liviu Aolaritei}, 
\author[First]{Saverio Bolognani},
\author[First]{Florian D\"orfler}

\address[First]{Automatic Control Laboratory, ETH Z\"urich, 
   Switzerland. \\ E-mail: aliviu@student.ethz.ch, \{bsaverio, dorfler\}@ethz.ch.}

\begin{abstract}                
We consider the problem of characterizing and assessing the voltage stability in power distribution networks.
Different from previous formulations, we consider the branch-flow parametrization of the power system state, which is particularly effective for radial networks.
Our approach to the voltage stability problem is based on a local, approximate, yet highly accurate characterization of the determinant of the power flow Jacobian.
Our determinant approximation allows us to construct a voltage stability index that can be computed in a fully scalable and distributed fashion.
We provide an upper bound on the approximation error, and we show how the proposed index outperforms other voltage indices that have been recently proposed in the literature.
\end{abstract}
\end{frontmatter}

\section{Introduction}

Operators of power distribution grids are facing unprecedented challenges caused by higher and intermittent consumers' demand, driven, among other things, by the penetration of electric mobility \citep{ClementNyns2010, Lopes2011}. Grid congestion is expected, as the demand gets closer to the hosting capacity of the network. 

One of the main phenomena that determines the finite power transfer capacity of a distribution grid is \emph{voltage instability} (see the recent discussion in \citealt{SimpsonPorco2016}).
The amount of power that can be transferred to the loads via a distribution feeder is inherently limited by the non-linear physics of the system.
In practice, as the grid load approaches this limit, increasingly lower voltages in the feeder are typically observed, followed by voltage collapse.

From the operational point of view, it is important to be able to identify operating conditions of the grid that are close to voltage collapse, in order to take the appropriate remedial actions. 
Although undervoltage conditions and voltage instability are related phenomena, it has been shown in \cite{Todescato2016} that it is not possible to identify the latter by simply looking at the feeder voltage levels.
Instead, many different indices have been proposed to quantify the distance of the grid from voltage collapse.
Most of them are based on the observation that the Jacobian of the power flow equations becomes singular at the steady state voltage stability limit (see the seminal work by \citealt{Tamura1988} and, even before, \citealt{Venikov1961}).
For a review of indices based on this approach, we refer to \cite{Chebbo1992} and to \cite{Gao1992}.

A geometric interpretation of the phenomena has been developed by \cite{Chiang1990}, and starting from 
\cite{Tamura1983} voltage collapse has been related to the appearance of bifurcations in the solutions of the nonlinear power flow equations.

More recently, semidefinite programming has been proposed as a tool to identify the region where voltage stability is guaranteed \citep{Dvijotham2015}.
The same region has been also characterized based on applications of fixed-point theorems (see \citealt{Bolognani2016} and references therein, and the extensions proposed in \citealt{Yu2015} and \citealt{Wang2016}).
Additionally, convex optimization tools have been used to determine sufficient condition for unsolvability (and thus voltage collapse) in \cite{Molzahn2013}.


All these works propose \emph{global indices}, in the sense that the knowledge of the entire system state is required at some central location, where the computation is performed.
Such a computation typically scales poorly with respect to the grid size, hindering the practical applicability of these methods.
Few exception include heuristic indices such as the one proposed in \cite{Vu1999}, which can be evaluated by each load based on local measurements.

The methodology that we propose in this paper builds on the aforementioned approach based on the singularity of the power flow Jacobian. Differently from other works, however, we adopt a branch flow model for the power flow equations \citep{Baran1989,Baran1989a,Farivar2013}.
This choice gives us a specific advantage, towards three results: first, we can reduce the  dimensionality of the problem via algebraic manipulation of the Jacobian of such equations; second, we can propose an approximation of the Jacobian-based voltage stability margin that is function of only the diagonal elements of the manipulated Jacobian, and is therefore computationally very tractable; finally, we can show how such an index can be computed in a completely distributed way, based on purely local measurements at the buses.
We derive an explicit bound for the approximation error, which is extremely small across the entire voltage stability region.
Based on that, we discuss how the proposed voltage stability index can be used in practice, and we show in numerical experiments how it outperforms other indices recently proposed in the literature.

The paper is structured in the following way.
In Section~\ref{sec:model} we recall the branch flow model, while in Section~\ref{sec:analysis} we explain how voltage stability can be assessed based on that model. 
In Section~\ref{sec:margin} we propose an approximate voltage stability index and we analyze the quality of the approximation. Finally, in Section~\ref{sec:applications}, we illustrate the result in simulations and we discuss the applicability of this approach to practical grid operation.

\section{Power distribution network model}
\label{sec:model}

Let $G = (N,E)$ be a directed tree representing a radial distribution network, where each node in $N = \{0,1,...,n\}$ represents a bus, and each edge in $E$ represents a line. Note that $|E|=n$. A directed edge in $E$ is denoted by $(i,j)$ and means that $i$ is the parent of $j$. For each node $i$, let $\delta(i) \subseteq N$ denote the set of all its children.
Node $0$ represents the root of the tree and corresponds to the distribution grid substation.
For each $i$ but the root $0$, let $\pi(i) \in N$ be its unique parent.

We now define the basic variables of interest.
For each $(i,j) \in E$ let $\ell_{ij}$ be the magnitude squared of the complex current from bus $i$ to bus $j$, and $s_{ij} = p_{ij} + \textbf{j}q_{ij}$ be the sending-end complex power from bus $i$ to bus $j$.
Let $z_{ij} = r_{ij} + \textbf{j}x_{ij}$ be the complex impedance on the line $(i,j)$.
For each node $i$, let $v_i$ be the magnitude squared of the complex voltage at bus $i$, and $s_i = p_i + \textbf{j}q_i$ be the net complex power demand (load minus generation) at bus $i$. 

Finally, we use the notation $\mathbf{1}$ and $\mathbf{0}$ for the vectors of all 1's and 0's, respectively.

\subsection{Relaxed branch flow model}

To model the power distribution network we use the relaxed branch flow equations proposed in \cite{Baran1989,Baran1989a,Farivar2013}%
\footnote{To make the model equations more compact, we adopted the convention 
$p_{\pi(0)0} = q_{\pi(0)0} = \ell_{\pi(0)0} = r_{\pi(0)0} = x_{\pi(0)0} = 0$.}
\begin{align*}
p_j &= p_{\pi(j)j} - r_{\pi(j)j}\ell_{\pi(j)j} - \sum\limits_{k \in \delta(j)}p_{jk}, \quad \forall j \in N\\
q_j &= q_{\pi(j)j} - x_{\pi(j)j}\ell_{\pi(j)j} - \sum\limits_{k \in \delta(j)}q_{jk}, \quad \forall j \in N\\
v_j &= v_i - 2(r_{ij}p_{ij} + x_{ij}q_{ij}) + (r_{ij}^2 + x_{ij}^2)\ell_{ij}, \  \forall (i,j) \in E \\
v_i \ell_{ij} &= p_{ij}^2 + q_{ij}^2, \quad \forall (i,j) \in E
\end{align*}

To write these equations in vector form, we first define the vectors $p$, $q$, and $v$, obtained by stacking the scalars $p_i$, $q_i$, and $v_i$, respectively, for $i \in N$.
Similarly we define $\overline{p}$, $\overline{q}$, $\ell$, $r$, and $x$, as the vectors obtained by stacking the scalars $p_{ij}$, $q_{ij}$, $\ell_{ij}$, $r_{ij}$, and $x_{ij}$, respectively, for $(i,j) \in E$.

In the following, we make use of the compact notation $[x]$, where $x \in \mathbb{R}^n$, to indicate the $n\times n$ matrix that has the elements of $x$ on the diagonal, and zeros everywhere else.

Finally, we define two $(0,1)$-matrices $A^i$ and $A^o$, where $A^i$ $\in$ $\mathbb{R}^{n+1 \times n}$ is the matrix which selects for each row $j$ the branch $(i,j)$, where $i = \pi(j)$, and $A^o \in \mathbb{R}^{n+1 \times n}$ is the matrix which selects for each row $i$ the branches $(i,j)$, where $j \in \delta(i)$.
Notice that $A := A^o-A^i$ is the incidence matrix of the graph.

The relaxed branch flow equations in vector form are:
\begin{equation}
\begin{split}
p &= A^i \big(\overline{p} - [r]\ell \big) - A^o \overline{p} \\
q &= A^i \big(\overline{q} - [x]\ell \big) - A^o \overline{q} \\
A^{i T} v &= A^{o T} v - 2 \big( [r]\overline{p} + [x]\overline{q} \big) + \big([r]^2 + [x]^2 \big) \ell \\
\left[A^{o T}v\right] \ell &= \left[ \overline{p}\right]\overline{p} + \left[\overline{q}\right]\overline{q}
\end{split}
\label{eq:branchflowmodelvector}
\end{equation}

We model node $0$ as a slack bus, in which $v_0$ is imposed ($v_0 = 1$ p.u.) and all the other nodes as PQ buses, in which the complex power demand (active and reactive powers) is imposed and does not depend on the bus voltage.
Therefore, the quantities $(v_0, p_{1\ldots n}, q_{1\ldots n})$ are to be interpreted as system parameters, and the relaxed branch flow model specifies $4n+2$ equations in $4n+2$ variables, $(\overline{p},\overline{q},\ell,v_{1\ldots n},p_0,q_0)$.

\section{Characterization of voltage stability}
\label{sec:analysis}

A \emph{loadability limit} of the power system is a critical operating point (as determined by the nodal power injections) of the grid, where the power transfer reaches a maximum value, after which the relaxed branch flow equations have no solution.
There are infinitely many loadability limits, corresponding to different demand configurations.
Ideally, the power system will operate far away from these points, with a sufficient safety margin.
On the other side, the \emph{flat voltage solution} (of the power flow equations) is the operating point of the grid where $v = \mathbf{1}$ and $p=q=\overline{p}=\overline{q}=\ell=\mathbf{0}$.
This point is voltage stable and the power system typically operates relatively close to it.

In the following, we recall and formalize the standard reasoning that allows to characterize loadability limits via conditions on the Jacobian of the power flow equations, and we specialize those results for the branch flow model that we have adopted.

\subsection{Jacobian of the power flow equations}

Based on the discussion at the end of Section~\ref{sec:model}, consider the two vectors
$$
u = 
\begin{bmatrix}
\overline{p} \\ \overline{q} \\ \ell \\ v_{1\ldots n} \\ p_0 \\ q_0
\end{bmatrix} \in \mathbb{R}^{4n+2}
\quad \text{and} \quad
\xi = 
\begin{bmatrix}
v_0 \\ p_{1\ldots n} \\ q_{1\ldots n}
\end{bmatrix} \in \mathbb{R}^{2n+1}
$$
corresponding to the system variables and the system parameters, respectively.
Then, the relaxed branch flow model \eqref{eq:branchflowmodelvector} can be expressed in an implicit form as
\begin{equation*}
	\varphi(u, \xi) = \mathbf{0}
\end{equation*}
From a mathematical point of view, a loadability limit corresponds to the maximum of a scalar function $\gamma(\xi)$ (to be interpreted as a measure of the total power transferred to the loads), constrained to the set $\varphi(u, \xi) = \mathbf{0}$ (the physical grid constraints).
\begin{align*}
	\max\limits_{u, \xi} \quad & \gamma(\xi)\\
	\text{subject to} \quad &  \varphi(u, \xi) = \mathbf{0}
\end{align*}

From direct application of the KKT optimality conditions, it results that in a loadability limit the power flow Jacobian $\varphi_u = \frac{\partial \varphi}{\partial u}$ becomes singular, i.e., $\det ( \varphi_u ) = 0$ (for details, see \citealt[Chapter 7]{Cutsem1998}). Based on this, we adopt the standard characterization for voltage stability of the grid, which we present in the following definition.

\emph{Definition.} (Voltage stability region). \ 
The voltage stability region of a power distribution network with one slack bus and $n$ PQ buses, described by the relaxed branch flow model, is the open region surrounding the flat voltage solution where the set of power flow solutions satisfy:
\begin{equation}
\det( \varphi_u ) \neq 0
\label{eq:detpositive}
\end{equation}

Although there might be other feasible regions, where the determinant of the power flow Jacobian is negative, the region characterized by \eqref{eq:detpositive} corresponds to the operating points of practical interest for the operation of the power system.

When the branch flow model is adopted, $\varphi_u$ takes the form
\begin{equation}
	\varphi_u = 
	\begin{bmatrix}
		-A & \mathbf{0}_{n+1 \times n} & -A^{i}[r] & \mathbf{0}_{n+1 \times n} & -\mathbf{e}_1 & \mathbf{0}_{n+1}\\
		\mathbf{0}_{n+1 \times n} & -A & -A^{i}[x] & \mathbf{0}_{n+1 \times n} & \mathbf{0}_{n+1} & -\mathbf{e}_1\\
		-2[r] & -2[x] & [r]^2 + [x]^2 & A_2^{T} & \mathbf{0}_n & \mathbf{0}_n\\
		2\left[\overline{p}\right] & 2\left[\overline{q}\right] & -\left[A^{o T}v\right] & -\left[\ell\right]A_2^{o T} & \mathbf{0}_n & \mathbf{0}_n
	\end{bmatrix}
\label{eq:pfj}
\end{equation}
where $A_2^{o}$ and $A_2$ are the matrices obtained by removing the first row from $A^{o}$ and $A$, respectively, and where $\mathbf{e}_1$ is the first canonical base vector.

Observe that the first three row blocks of $\varphi_u$ are constant, while the last row block depends linearly on the variables $\overline{p}, \overline{q}, v$ and $\ell$.

\subsection{Reduced power flow Jacobian}


We define the following $n \times n$ matrix, that we denote as the \emph{reduced power flow Jacobian}.
\begin{multline}
	 \varphi_u^{'} = \left[ A^{o T}v \right] + 2\left[\overline{p}\right]A_2^{-1}[r] + 2\left[\overline{q}\right]A_2^{-1}[x]  \\
- [\ell] A_2^{o T}(A_2^T)^{-1}  \left( [r]^2 + 2 [r] A_2^{-1} [r] + [x]^2 + 2 [x] A_2^{-1} [x] \right)
\label{eq:rpfj}
\end{multline}

In the following, we provide a key theorem that shows the merits of the reduced power flow Jacobian.

\begin{thm}
Consider the power flow Jacobian \eqref{eq:pfj} and the reduced power flow Jacobian \eqref{eq:rpfj}
of a power distribution network with one slack bus and $n$ PQ buses, described by the relaxed branch flow model.
The following statements hold.
\begin{itemize}
\item[i)] $\det( \varphi_u ) =  \det( \varphi_u' )$.
\item[ii)] $\det( \varphi_u' ) > 0$ in the voltage stability region.
\item[iii)] $\det( \varphi_u') = 0 \ \Leftrightarrow \ \exists \ i \in \delta(0): \det( \varphi_{u,i}' ) = 0$, 
where $ \varphi_{u,i}'$ is the reduced power flow Jacobian of the tree composed by node $0$ and the subtree rooted by child $i$ of node $0$.
\end{itemize}
\label{thm:voltagestabilityregion}
\end{thm}
\begin{pf}
For \emph{i) } and \emph{iii) } only a sketch of the proof is provided. The full details are available in \cite{Aolaritei2016}. \\
\emph{i)\ }
Observe that the last two columns of $\varphi_u$ are the canonical vectors $\mathbf{e}_1$ and $\mathbf{e}_{n+2}$ of $\mathbb{R}^{4n+2}$.
Thus, if we eliminate these columns together with the $1^{st}$ and $n+2^{nd}$ rows we obtain a new matrix, $\varphi_u^{*}$ of dimensions $4n \times 4n$ whose determinant is equal to $(-1)^n \det(\varphi_u)$.
We next prove that $\det ( \varphi_u^{'} ) =  (-1)^n \det(\varphi_u^{*}$ ).
To do so, we apply  Schur complement twice on the matrix $\varphi_u^{*}$.
After some very basic matrix manipulation the result is obtained.\\
\emph{ii)\ } 
In the flat voltage solution we have that $\varphi_u^{'} = \left[ A^{o T}v \right]= \left[ A^{o T} \mathbf{1} \right] = I$,
and therefore $\det(\varphi_u^{'}) = 1$.
Thus there exists a solution in the voltage stability region where the determinant of the power flow Jacobian is positive.
Moreover, we know that in a loadability limit, $\det(\varphi_u^{'}) = 0$, and that the determinant is a continous function of the grid variables.
Therefore, in order to remain in the voltage stability region, the determinant needs to remain positive.\\
\emph{iii)\ }
By re-indexing the nodes of the network, $\varphi_u^{'}$ can be transformed in a block diagonal matrix, where each block depends only on node 0 and the subtree rooted by one child of node $0$.
$\hfill \square$
\end{pf}

Theorem~\ref{thm:voltagestabilityregion} shows that the reduced power flow Jacobian $\varphi_u'$ is an effective tool for the characterization of the voltage stability region, and for the voltage stability analysis of a distribution grid. 
In particular, i) shows that studying the reduced power flow Jacobian is completely equivalent to studying the original power flow Jacobian, when we are interested in its singularity. 
ii) provides a more precise characterization of the region where the grid voltages are stable.
Finally, iii) explains how the dimensionality of the problem of computing the determinant of the power flow Jacobian can be further reduced, if the root (node $0$) has more than one child.



\section{Voltage stability analysis}
\label{sec:margin}

In this section we first propose an approximation of the determinant of the reduced power flow Jacobian that is amenable to scalable and distributed computation, when measurements of the grid variables are available. Then, based on this approximation, we propose a voltage stability index to quantify the distance of the power system from voltage collapse.

\subsection{Mathematical preliminaries on matrix theory}

Given $A \in \mathbb{R}^{n \times n}$, we denote by $A_\text{diag}$ and $A_\text{off}$ the matrices that contain only the diagonal and off-diagonal elements of $A$, respectively. We denote by $\rho = \rho(A)$ its spectral radius, i.e. the maximum norm of its eigenvalues.

\emph{Definition.} A matrix $A \in \mathbb{R}^{n \times n}$ is a $Z$-matrix if $A = \alpha I - B$, where $\alpha$ is a real number and $B$ is a nonnegative matrix. The set of all $n \times n$ $Z$-matrices is denoted by $Z_{<n>}$.


\emph{Definition.} A matrix $A \in \mathbb{R}^{n \times n}$ is an $\omega$-matrix if:
\begin{enumerate}
	\item Each principal submatrix of A has at least one real eigenvalue.
	\item If $S_1$ is a principal submatrix of $A$ and $S_{11}$ a principal submatrix of $S_1$ then $\lambda_{min}(S_1) \leq \lambda_{min}(S_{11})$, where $\lambda_{min}$ denotes the smallest real eigenvalue.
\end{enumerate}
The set of all $n \times n$ $\omega$-matrices is denoted by $\omega_{<n>}$.

\emph{Definition.} A matrix $A \in \mathbb{R}^{n \times n}$ is a $\tau$-matrix if it is an $\omega$-matrix and $\lambda_{min}(A) \geq 0$. The set of all $n \times n$ $\tau$-matrices is denoted by $\tau_{<n>}$.


The following technical results will be used.

\begin{thm}[\citealt{Mehrmann1984}]
$Z_{<n>}$ $\subseteq$ $\omega_{<n>}$.
\label{thm:zw}
\end{thm}

\begin{thm}[\citealt{Engel1975}]
If $A \in \mathbb{R}^{n \times n}$ is a $\tau$-matrix then:
\begin{equation*}
	\det(A) \leq \det(A_\text{diag})
\end{equation*}
\label{thm:hadamard}
\end{thm}

\begin{thm}[\citealt{Ipsen2011}]
Given $A \in \mathbb{R}^{n \times n}$, if $A_\text{diag}$ is nonsingular and $\rho = \rho(A_\text{diag}^{-1}A_\text{off}) < 1$, then:
\begin{equation}
	\left|\ln(\det(A)) - \ln(\det(A_\text{diag}))\right| \ \leq \ - \rho n \, \ln(1 - \rho)
\label{eq:detapprox}
\end{equation}
\label{thm:detapprox}
\end{thm}

\subsection{Determinant approximation}

Direct inspection of the reduced power flow Jacobian $\varphi_u^{'}$ shows that, for realistic parameter values and operating conditions, its off-diagonal elements (and in particular its lower-diagonal elements) are significantly smaller than the diagonal elements.
The approximation proposed in this paper consists in ignoring them, and requires the following assumption.

\begin{assum}
All PQ buses in the network have positive active and reactive power demand.
\label{ass:loads}
\end{assum}

This assumption ensures that $p_{ij}, q_{ij} \geq 0$  $\forall (i,j) \in E$, although it is not a necessary condition for that to hold.
In practical terms, having positive power demands everywhere corresponds to the most unfavorable case for voltage stability, and there is little loss of generality in assuming that in this analysis.
Based on this assumption, in the remaining of this paper we will refer to the nodes $1,...,n$ as \emph{PQ loads}.

\begin{figure}[tb]
	\begin{center}
		\includegraphics[width=0.48\columnwidth]{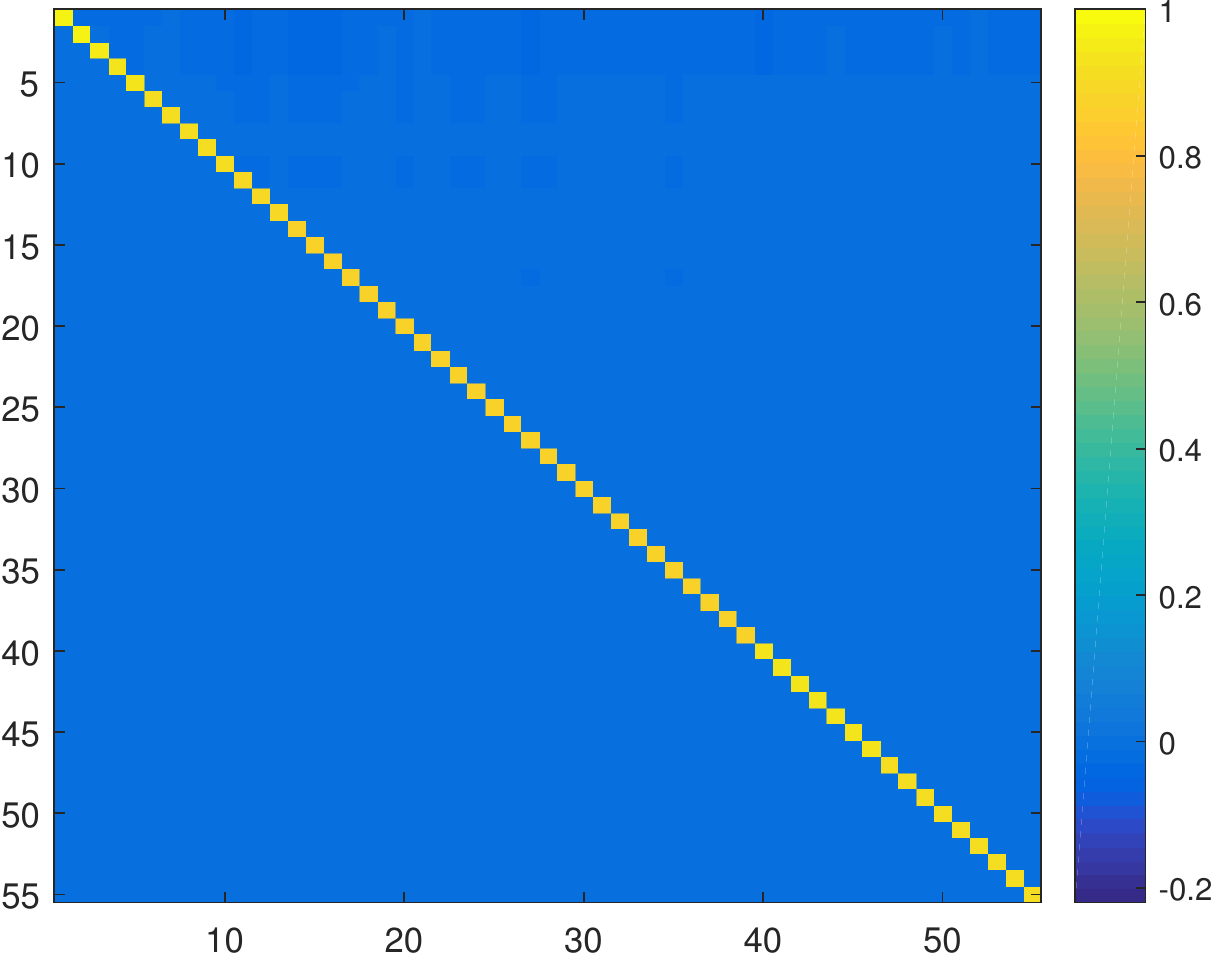}\hspace{\stretch{1}}
		\includegraphics[width=0.48\columnwidth]{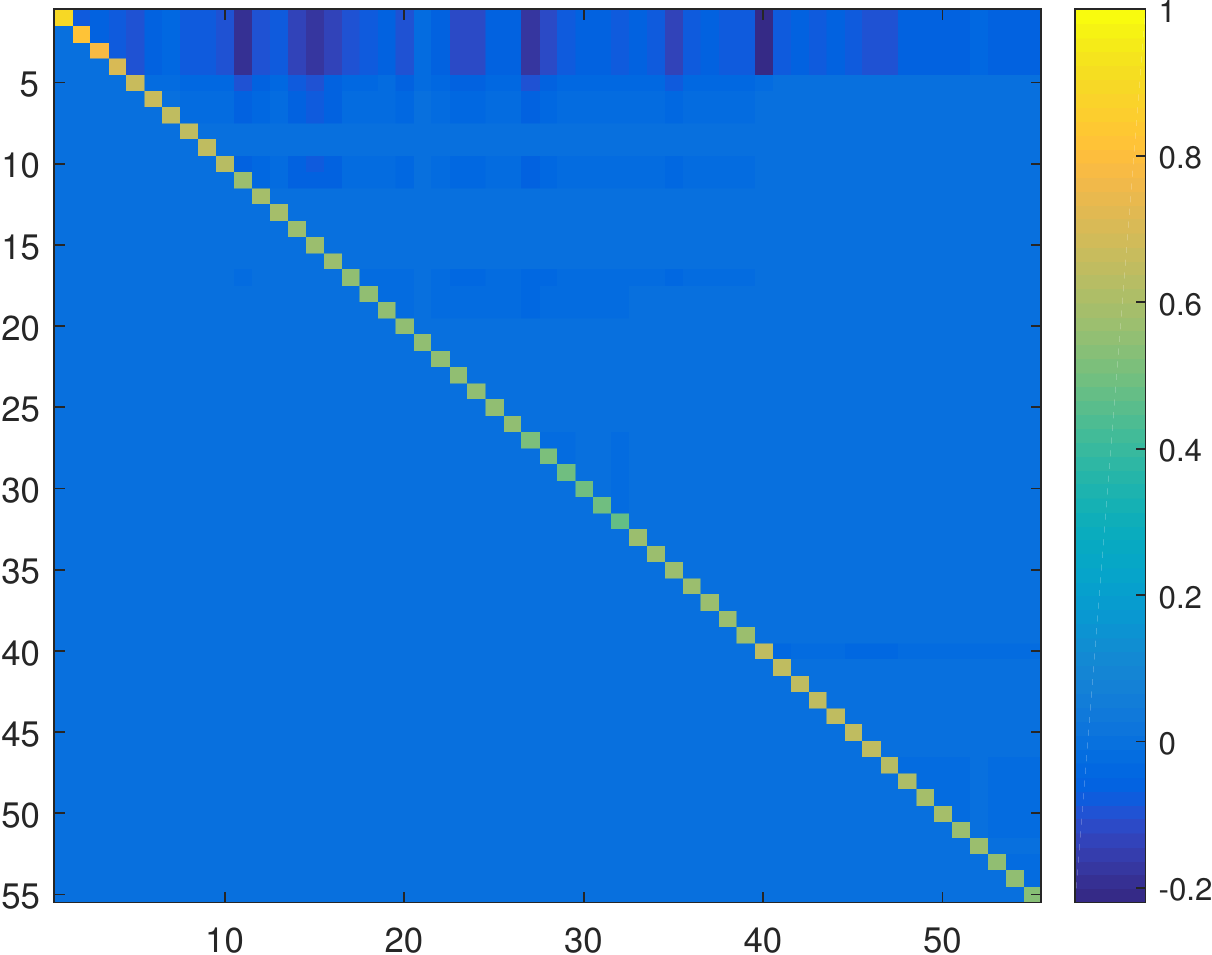}
	\caption{Data in the reduced power flow Jacobian}
	\label{fig:datajacobian}
	\end{center}
\end{figure}

In Fig.~\ref{fig:datajacobian} we represent the numerical value of $\varphi_u^{'}$ for two levels of loadability of a 56-bus distribution grid (described in detail in Section~\ref{sec:applications}).
In the left panel, the operating point of the system is close to the flat voltage solution, while in the right panel, the grid is operated close to a loadability limit.

The diagonal elements of $\varphi_u^{'}$ are equal to 
\begin{equation}
\varphi_{u,jj}' = v_i - 2p_{ij}r_{ij} - 2q_{ij}x_{ij} - 2\ell_{ij}(r_{ij}\overline{r_{0i}} + x_{ij}\overline{x_{0i}})
\label{eq:diagelems}
\end{equation}
where $i = \pi(j)$ and $\overline{r_{0i}}$ is the sum of the resistances of the lines connecting node $0$ to node $i$ (and similarly for $\overline{x_{0i}}$).


By ignoring the off-diagonal elements, an approximation of $\det(\varphi_u^{'})$ is obtained as the product of the elements on the diagonal defined in \eqref{eq:diagelems}:
\begin{equation}
	{\det}_\text{approx} = \prod_{(i,j)\in E} \varphi_{u,jj}'
	\label{eq:detapproxprod}
\end{equation} 

In the next Lemma, we prove that the approximation is an upper bound for the true determinant.

\begin{lem} For a power distribution network with one slack bus and $n$ PQ loads described by the relaxed branch flow model, in the voltage stability region the determinant of the reduced power flow Jacobian satisfies
\begin{equation*}
	0 < \det(\varphi_u^{'}) \leq {\det}_\text{approx}
\end{equation*}
\label{lem:lowerbound}
\end{lem}
\begin{pf}
Having $p_{ij}, q_{ij} \geq 0 \ \forall (i,j) \in E$ ensures that the off-diagonal elements of $\varphi_u^{'}$ are nonpositive.
Thus, $\varphi_u^{'}$ is a $Z$-matrix and, from Theorem~\ref{thm:zw}, $\varphi_u^{'}$ is also an $\omega$-matrix.
Recall that an $\omega$-matrix is a nonsingular $\tau$-matrix if and only if its smallest real eigenvalue is positive.
But this is exactly what we require for voltage stability.
To see this, notice that in the flat voltage solution all the eigenvalues of $\varphi_u^{'}$ are real and equal to 1.
Thus the determinant becomes zero for the first time when the smallest real eigenvalue becomes zero.
Notice that there is always at least one real eigenvalue, since $\varphi_u^{'}$ is an $\omega$-matrix.
Therefore, in the voltage stability region, $\varphi_u^{'}$ is a $\tau$-matrix.
The result follows from Theorem~\ref{thm:hadamard}.
$\hfill \square$
\end{pf}

Numerical experiments show that the approximation is exact only in the flat voltage solution, though the approximation error is almost negligible (see Section 5).
With positive power demands, $\det(\varphi_u^{'}) < \det_\text{approx}$.



\subsection{Voltage stability index}
\label{ssec:vsi}

Based on Theorem~\ref{thm:voltagestabilityregion}, the voltage stability region is defined as the region where $\det(\varphi_u^{'}) > 0$.
In practical terms, the grid operator has to identify a threshold $\beta>0$ and impose that $\det(\varphi_u^{'}) \geq \beta$ as a practical voltage stability measure.
In order to make full use of the capacity of the grid, the value $\beta$ needs to be chosen such that, when $\det(\varphi_u^{'}) = \beta$, the operating point of the grid is very close to a loadability limit. 
From numerical experiments, it is evident that a proper choice of $\beta$ intrinsically depends on the size of the network.
To gain some intuition about this, recall that the determinant of a matrix is equal to the product of its eigenvalues. 
In the flat voltage solution, all the eigenvalues of $\varphi_u^{'}$ are equal to 1.
As soon as the power demands increase, the eigenvalues start moving towards the origin.
Since the number of eigenvalues is equal to the size of $\varphi_u^{'}$, and thus to the size of the network, it is clear that bigger networks are associated to exponentially smaller determinants. 

Based on this intuition, we propose 
\begin{equation*}
\text{VSI} := \frac {\ln(\det(\varphi_u^{'}))}{n}
\end{equation*}
as a \emph{voltage stability index}.
Thus, for some threshold $\beta > 0$, the practical voltage stability measure becomes
\begin{equation}
	\text{VSI} \geq \frac{\ln(\beta)}{n} := \text{VSI}_\text{min}
\label{eq:vsimin}
\end{equation}

Following the determinant approximation proposed in \eqref{eq:detapproxprod}, we then define the voltage stability index approximation
\begin{equation*}
\text{VSIA} := \frac {\ln\left(\det_\text{approx}\right)}{n}
\end{equation*}

In the following remark we point out an interesting and useful property of this voltage stability index approximation.

\begin{rem}[Distributed computation of the VSIA]
\ \\Notice that $\varphi_{u,jj}'$ is only function of the local state variables relative to the edge $(i,j)$, where $i=\pi(j)$.
More precisely, $\varphi_{u,jj}'$ can be computed in a distributed way from measurements performed at bus $i$ and on the power lines that leave the same bus: $v_i$, $p_{ij}$, $q_{ij}$ and $\ell_{ij}$.
Once each node $i$ has computed $\varphi_{u,jj}'$ for each children $j \in \delta(i)$, the computation of the $\text{VSIA}$ amounts to simply evaluating the arithmetic mean of the terms $\ln\left(\varphi_{u,jj}'\right)$ for all $(i,j) \in E$.
The arithmetic mean of these nodal quantities can then be computed via scalable fully distributed algorithms such as consensus algorithms \citep{Olfati-Saber2004}.
\end{rem}

\subsection{Approximation error}

In this section we study the approximation error between $\text{VSI}$ and $\text{VSIA}$.
To do so, we need the following lemma.

\begin{lem} 
In a power distribution network with one slack bus and $n$ PQ loads described by the relaxed branch flow model, in the voltage stability region,  the reduced power flow Jacobian satisfies the following:
\begin{itemize}
\item[i)] $\varphi_{u,\text{diag}}'$ is positive definite
\item[ii)] $\rho(\varphi_{u,\text{diag}}^{' -1}\varphi_{u,\text{off}}^{'}) < 1$
\label{lem:rho}
\end{itemize} 
\end{lem}
\begin{pf}
\emph{i)\ }
The two facts, $\varphi_{u,\text{diag}}' = I$ in the flat voltage solution and $\det(\varphi_{u,\text{diag}}^{'}) > 0$ in the voltage stability region (via Lemma~\ref{lem:lowerbound}), ensure that the elements on the diagonal remain positive.\\
\emph{ii)\ }
Since $\varphi_{u}^{'} = \varphi_{u,\text{diag}}^{'}(I + \varphi_{u,\text{diag}}^{' -1}\varphi_{u,\text{off}}^{'})$, we have that $\det(\varphi_{u}^{'}) = \det(\varphi_{u,\text{diag}}^{'})\det(I + \varphi_{u,\text{diag}}^{' -1}\varphi_{u,\text{off}}^{'})$.
In the flat voltage solution, $\varphi_{u,\text{diag}}^{' -1}\varphi_{u,\text{off}}^{'} = \mathbf{0}_{n \times n}$ and in a loadability limit, $\det(I + \varphi_{u,\text{diag}}^{' -1}\varphi_{u,\text{off}}^{'}) = 0$. Thus, the power grid becomes unstable when an eigenvalue of $\varphi_{u,\text{diag}}^{' -1}\varphi_{u,\text{off}}^{'}$ arrives at $-1$.
Now, since $-\varphi_{u,\text{diag}}^{' -1}\varphi_{u,\text{off}}^{'}$ is non-negative, it has a positive real eigenvalue equal to the spectral radius $\rho(-\varphi_{u,\text{diag}}^{' -1}\varphi_{u,\text{off}}^{'})$ (Perron-Frobenius Theorem). 
Therefore, $\varphi_{u,\text{diag}}^{' -1}\varphi_{u,\text{off}}^{'}$ has a negative real eigenvalue with magnitude equal to $\rho(\varphi_{u,\text{diag}}^{' -1}\varphi_{u,\text{off}}^{'})$. Hence, this is the eigenvalue that first arrives in $-1$. This implies that in the voltage stability region, $\rho(\varphi_{u,\text{diag}}^{' -1}\varphi_{u,\text{off}}^{'}) < 1$.
$\hfill \square$
\end{pf}

In the following Lemma we give an exact expression for the approximation error.

\begin{lem}In a power distribution network with one slack bus and $n$ PQ loads described by the relaxed branch flow model, in the voltage stability region we have:
\begin{equation}
	\text{VSIA} - \text{VSI} = \frac{\trace(\sum_{i=2}^{\infty}\frac{(-\varphi_{u,\text{diag}}^{' -1}\varphi_{u,\text{off}}^{'})^i}{i})}{n}
\end{equation}
\label{eq:traceformula}
\end{lem}
\begin{pf}
We have that $\ln(\det(\varphi_{u}^{'})) = \ln(\det(\varphi_{u,\text{diag}}^{'})) + \ln(\det(I + \varphi_{u,\text{diag}}^{' -1}\varphi_{u,\text{off}}^{'}))$.
As $\rho(\varphi_{u,\text{diag}}^{' -1}\varphi_{u,\text{off}}^{'}) < 1$ we know that $\ln(\det(I + \varphi_{u,\text{diag}}^{' -1}\varphi_{u,\text{off}}^{'})) = \trace(\ln(I + \varphi_{u,\text{diag}}^{' -1}\varphi_{u,\text{off}}^{'})) = -\trace(\sum_{i=1}^{\infty}(-\varphi_{u,\text{diag}}^{' -1}\varphi_{u,\text{off}}^{'})^i/i)$.
To conclude, notice that $\trace(\varphi_{u,\text{diag}}^{' -1}\varphi_{u,\text{off}}^{'}) = 0$.
$\hfill \square$
\end{pf}


In Section 5 we show that the approximation is almost exact in the voltage stability region.
Since the terms in the above sum are all positive, they are very small and they decay quickly to zero.
The numerical value of the right hand side of \eqref{eq:traceformula} has been plotted in Fig.~\ref{fig:errorbounds}, for a the test distribution feeder described in Section~\ref{sec:applications}, and for different load levels.

In the following theorem, we present the main result on the quality of the proposed voltage stability index approximation.

\begin{thm}
In a power distribution network with one slack bus and $n$ PQ loads described by the relaxed branch flow model, in the voltage stability region we have:
\begin{equation}
\text{VSI} \leq \text{VSIA} \leq \text{VSI}  -\rho \, \ln(1 - \rho)
\label{eq:mainresult}
\end{equation}
where $\rho = \rho(\varphi_{u,\text{diag}}^{' -1}\varphi_{u,\text{off}}^{'})$.
\label{thm:mainresult}
\end{thm}
\begin{pf}
The first inequality descends from Lemma~\ref{lem:lowerbound}.
The second inequality is proved by applying Theorem~\ref{thm:detapprox}, using what we proved in Lemma~\ref{lem:rho}.
$\hfill \square$
\end{pf}

We conclude this section by presenting the following conjecture.

\begin{conj}
In~\cite{Ipsen2011}, the authors illustrate that the pessimistic factor in the approximation bound of Theorem~\ref{thm:detapprox} is given by the factor $n$ that appears in \eqref{eq:detapprox}.
They found that replacing $n$ by the number of eigenvalues whose magnitude is close to the spectral radius makes the bound tight.
In our simulations we found that there is generally only one eigenvalue with magnitude close to the spectral radius.
This would imply that the result that we presented in Theorem~\ref{thm:mainresult} can be tightened to
\begin{equation}
 \text{VSI} \leq \text{VSIA} \leq \text{VSI}  - \frac{1}{n}\rho \, \ln(1 - \rho)
\label{eq:conjecture}
\end{equation}
\label{conjecture}
\end{conj}

This tighter bound on the approximation error has always revealed to be true in our simulations, as illustrated in the next section.

\section{Numerical validation and comparison}
\label{sec:applications}

\subsection{Numerical validation of the VSI approximation}

In this section we assess the quality of the proposed voltage stability index approximation via numerical simulations.
We consider a 56-bus distribution network, obtained from the three-phase backbone of the IEEE123 test feeder.
The details of the testbed are available in \cite{github_approx-pf}. Power flow equations have been solved via MatPower \citep{Zimmerman2011}.

\begin{figure}[t]
	\centering
		\includegraphics[width=0.95\columnwidth]{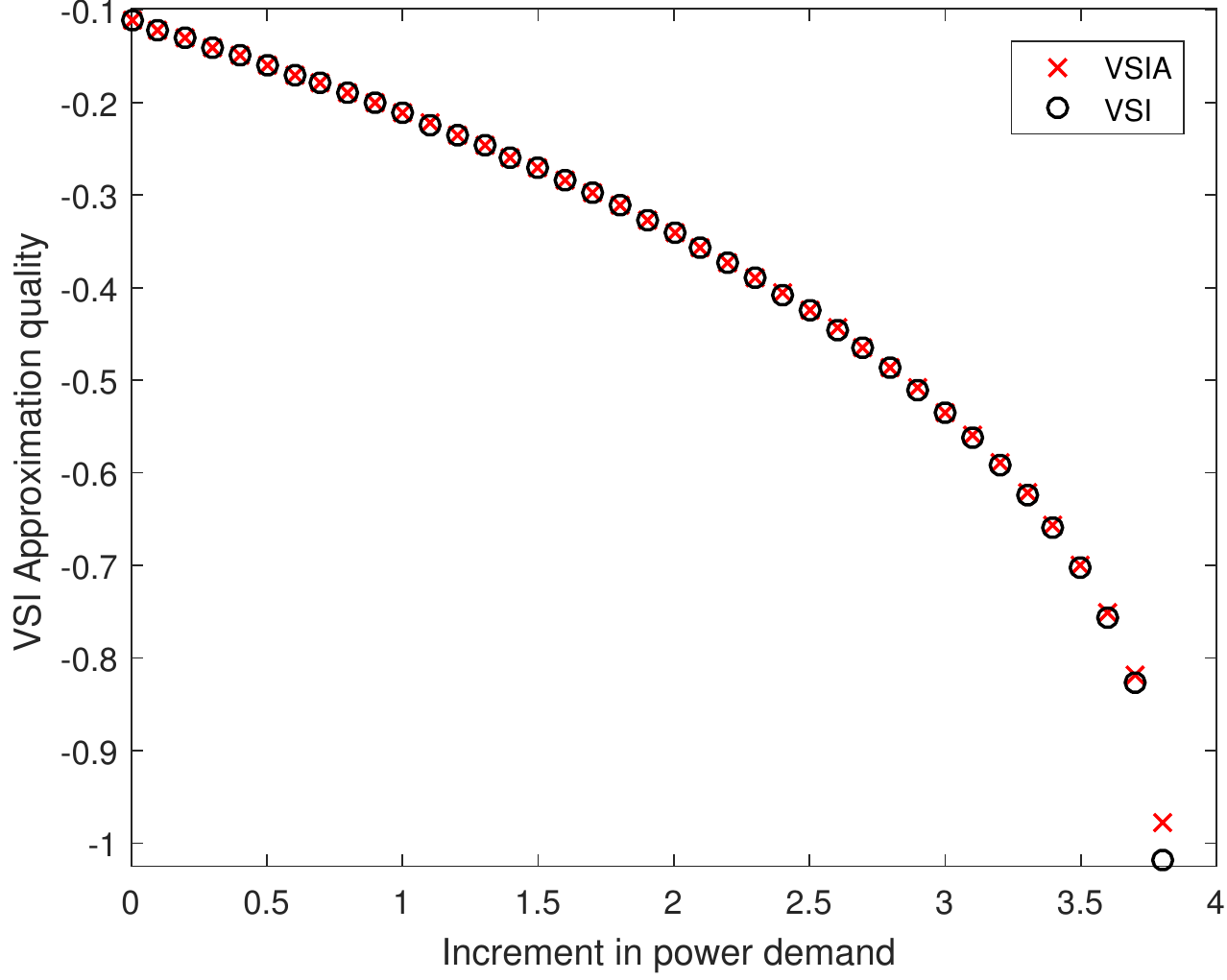}
	\caption{Voltage stability index and its approximation, for a series of increasing demand levels, from flat voltage to voltage collapse.}
	\label{fig:vsia}
\end{figure}

In Fig.~\ref{fig:vsia} we represent the voltage stability index (VSI) and the voltage stability index approximation (VSIA) when the system is operated at a series of increasing power demands.
We start from an operating point very close to the flat voltage solution, and we increase the active and reactive power demand at four different buses in the grid until the Jacobian becomes singular and the Newton's method employed for the solution of the power flow equations cannot proceed.
Observe, that the proposed VSI approximation is almost exact up to very close to the loadability limit.

\begin{figure}[t]
	\centering
		\includegraphics[width=0.95\columnwidth]{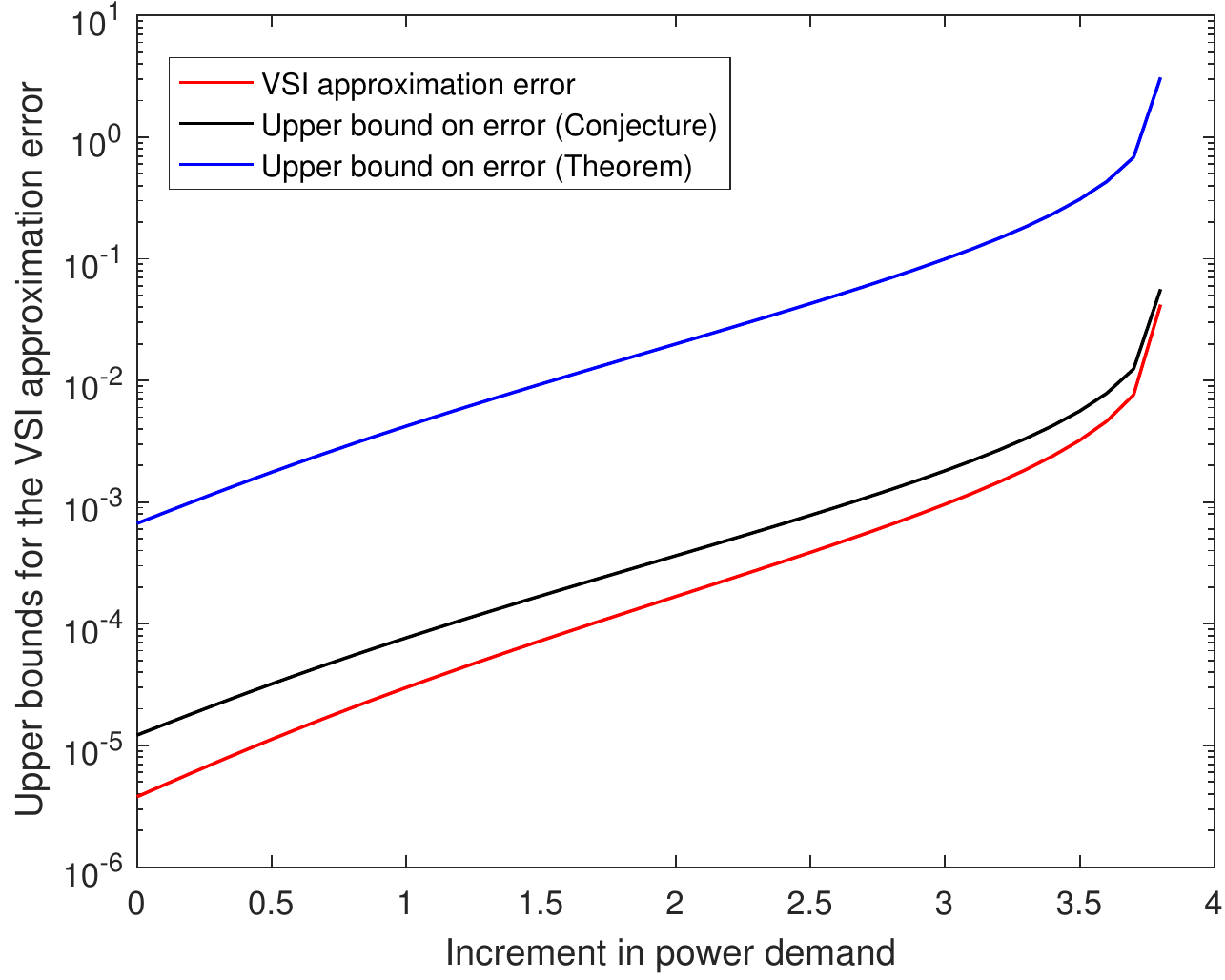}
	\caption{Comparison of the voltage stability index approximation errors with the two proposed error bounds.}
	\label{fig:errorbounds}
\end{figure}

In Fig.~\ref{fig:errorbounds} we represent the VSI approximation error, together with the bounds presented in Theorem~\ref{thm:mainresult} and Conjecture~\ref{conjecture}.
Observe that the approximation error is quite small in either case, and it follows the conjectured bound \eqref{eq:conjecture} rather than the bound \eqref{eq:mainresult}.

More simulations can be found in \cite{Aolaritei2016}, and show how the quality of the approximation is consistently good across different power demands configurations.

\subsection{Comparison of practical voltage stability indices}

Recall from Section~\ref{ssec:vsi} that we propose $\text{VSI} \ge \text{VSI}_\text{min}$ as a voltage stability measure, where $\text{VSI}_\text{min}$ has to be decided in order to characterize an operating condition close to the loadability limit of the grid.
It can be seen in Fig.~\ref{fig:vsia} that when $\text{VSI}=-1$, its negative slope is already extremely steep, meaning that for a very small increase in power demand the system would become unstable. 
Preliminary numerical investigation has shown that this threshold for $\text{VSI}$ is valid for a diverse range of grid sizes and topologies.
Notice that such a limit corresponds to an exponentially decreasing threshold for the determinant of the power flow Jacobian, i.e., $\det(\varphi_u') \ge e^{-n}$.

In practical terms, however, when the VSI is to be used as a tool for the assessment of the distance from voltage collapse, a more conservative value of $\text{VSI}_\text{min}$ is to be chosen. 
In the following, we choose a slightly more conservative limit ($\text{VSI}_\text{min} = -0.8$) in order to present a comparison between the proposed VSI and three other indices that have been recently proposed in the literature.
Observe from Fig.~\ref{fig:vsia} that the approximation is extremely precise when $\text{VSI}$ is larger than $-0.8$. 
This suggests that the VSIA can be safely used instead, enabling a fast, scalable, and distributed assessment of the voltage stability of the grid.

The first two indices that we consider have been  proposed in \cite{Bolognani2016} and in \cite{SimpsonPorco2016}, and they involve open-circuit load voltages, the grid impedance matrix (or a specific norm of it), and nodal power injections.
The third index has been presented in \cite{Wangtoappear}, and it requires the knowledge of the impedance matrix of the grid and of phasorial measurements of the bus voltages.

For each criterion we evaluated the proposed voltage stability index in a low-load operating point (very close to the flat voltage profile) and in the operating point in which our VSIA becomes equal to $-0.8$ (corresponding to what we defined as the practical voltage stability limit). 
Since the method proposed in \cite{SimpsonPorco2016} is based on the decoupled reactive power flow equations, for the comparison with their method we used only reactive power demands.


We obtained the following values, showing how the proposed voltage stability index approximation is in fact an effective tool for the precise assessment of the distance of the system from voltage collapse. 
The other indices reach their threshold value before our index does, showing that they are more conservative, and therefore result in a less efficient use of the given distribution grid.

\begin{center}
\footnotesize
\begin{tabular}{@{}lcc@{}}
\toprule
Criterion & Low load & $\text{VSIA} = -0.8$\\
\midrule
$\text{VSIA} > -0.8$ & -0.10 & -0.80\\
$\text{VSI$^{\text{Bolognani}}$} > 0$ & 0.51 & -4.58 \\
$\text{VSI$^{\text{Simpson}}$} < 1$ & 0.32 & 1.34 \\
$\text{VSI$^{\text{Wang}}$} > 1$ & 10.81 & 0.97 \\
\bottomrule 
\end{tabular}
\end{center}

\section{Conclusions}

In this paper we have presented a voltage stability index for power distribution networks, for which an accurate approximation is available. 
Bounds on the quality of this approximation have been mathematically derived, and the accuracy have been validated in simulations.
Notably, the approximate voltage stability index can be computed in a scalable and distributed way by agents that can measure local variables at each bus.
Based on this observation, we envision three possible applications for which the proposed approach can bring a significant contribution.
\begin{itemize}
\item As an online voltage stability monitoring tool, when the necessary quantities are measurable at the buses, the VSIA can be computed asynchronously via standard tools from multi-agent average consensus.
\item In optimal power flow programming, whenever the problem is expressed via the branch flow model, the VSIA can be used as a computationally efficient barrier function to maintain the solution of the problem inside the region of voltage stability.
\item In numerical algorithm for the construction of power flow feasibility sets that are based on the nonsingularity of the power flow Jacobian (as in \citealt{Dvijotham2015}), the proposed approximation can be used to avoid expensive determinant computations and improve scalability to larger networks.
\end{itemize}
Preliminary numerical investigation shows that the voltage stability index approximation remains extremely accurate even in the presence of generators (i.e., positive power injections, which violate Assumption~\ref{ass:loads}).
An extension to this more general case is currently under development, together with a numerical assessment of the effectiveness of the proposed index on various distribution test feeders.

\bibliography{voltagestability} 

\end{document}